\documentclass{ifacconf}

\usepackage{amssymb}
\usepackage{bm}
\usepackage{amsmath}
\usepackage{graphicx}
\usepackage{xcolor}
\usepackage{mathtools}
\usepackage{natbib}

\newcommand{\J}{{\cal J}}

\begin{document}
\begin{frontmatter}

\title{Network Inference using Sinusoidal Probing}

\thanks[footnoteinfo]{This work has been supported by the Swiss National Science Foundation under grants 200020\_182050 and P400P2\_194359.}

\author{Robin Delabays$^{1,2}$ and Melvyn Tyloo$^{3,4}$} 

\address{$^1$ Automatic Control Laboratory, Swiss Federal Institute of Technology, CH-8092 Zurich, Switzerland. \\
$^2$ Center for Control, Dynamical Systems and Computation, UCSB, Santa Barbara, CA 93106-5070 USA. \\
$^3$ School of Engineering, University of Applied Sciences of Western Switzerland HES-SO, CH-1951 Sion, Switzerland. \\
$^4$ Institute of Physics, \'Ecole Polytechnique F\'ed\'erale de Lausanne (EPFL), CH-1015 Lausanne, Switzerland.}

\begin{abstract}
 The aim of this manuscript is to present a non-invasive method to recover the network structure of a dynamical system. 
 We propose to use a controlled probing input and to measure the response of the network, in the spirit of what is done to determine oscillation modes in large electrical networks. 
 For a large class of dynamical systems, we show that this approach is analytically tractable and we confirm our findings by numerical simulations of networks of Kuramoto oscillators. 
 Our approach also allows us to determine the number of agents in the network by probing and measuring a single one of them. 
\end{abstract}


\end{frontmatter}

%===============================================================================

\section{Introduction}
Complex networks are the medium for interactions in many natural and man-made systems. 
Such realizations range from the scale of people exchanging opinions on a social networks and power transmission on electrical grids to interacting molecules in chemical reactions and pacemaker cells~[\cite{Boc06,New18}]. 
The way individual elements are coupled together primarily impacts the overall dynamics of network-coupled systems. 
However, in many cases, characteristics of the interaction network are not known exactly, or even not known at all. 
Inference techniques to uncover coupling between individual units and even units with themselves are therefore highly desirable. 
We distinguish mainly two types of such methods, 
(i) one can observe a dynamical system subjected to uncontrolled operational condition; or  
(ii) one can directly disturb the system and observe its reaction, which is what is done in this manuscript. 
In real applications, one has to be careful on the nature of the method. 
Introducing a disturbance into the system can have dramatic effect on its operation. 
For instance, \cite{Fur19} propose a method based on resonance of the network when subject to periodic disturbance, and \cite{Tim07} suggests to drive the system away from its operating state. 
Such (potentially invasive) methods could alter the operational state of the system. 

To avoid a strong impact on the system, one can adopt the strategy of simply observing the system in its normal operation. 
This type of approach lead to successful and elegant results, e.g., leveraging the response of the system to external noise [\cite{Ren10,Tyl20b}]. 
But as there is no free lunch, these methods come with assumptions on the noise characteristics that are not necessarily met in general, namely on correlation time and uniformity over the system. 
Observing the relaxation of a system to its steady state, \cite{Mau17} extract its spectral moments, but do not directly reconstruct the network. 
Furthermore, their method performs exactly for linear dynamics, but relies on Dynamical Mode Decomposition for nonlinear systems, which might require a very large number of measurements for an accurate estimation. 

Half way between intrusive disturbances and passive observation, we will consider the injection of a small probing signal, which we will qualify as "non-invasive", in order to leave the operational state as unaffected as possible in the spirit of \cite{Pie10}. 
To this day, most of such techniques, applied to large power grids, relied of very few measurement points and aimed only at identifying resonance modes of the network, and not the whole network structure [\cite{Dos13}].
In this manuscript, we propose such a non-invasive inference technique. 
It relies on rather mild assumptions on the nature of the interaction between agents of the network, and does not need to know it, on the contrary to other methods [e.g., \cite{Yu06}]. 
Our method applies generically to any network and does not need any knowledge of its charateristics, as for instance in~\cite{Yeu02}. 
By adding a controlled sinusoidal input at single nodes, referred to as \emph{probing signal}, we are able to reconstruct the interaction network by measuring the response to the probing at the other nodes of the network. 
The same approach allows us to determine the number of nodes in the network by probing and measuring the response of the network at a single node. 
This improves significantly on previous measurement-based methods to determine the number of nodes in a network [\cite{Hae19}].
In contrast to more exhaustive works [e.g., \cite{Mat10,Dan12} and other publication of the same authors], where the authors aim at identifying both network structure, transfer functions, and internal node dynamics, our method focusses on the interaction network only, avoiding to rely on an error minimization, which cannot always be computed in closed form.

\section{Preliminaries}
Let us consider a general network of $n$ coupled agents with first-order dynamics and diffusive coupling
\begin{align}\label{eq:dyn}
 \dot{y}_i &= \omega_i - \sum_j a_{ij}f_{ij}(y_i-y_j) + \xi_i\, , & i &= 1,...,n\, ,
\end{align}
where $y_i\in \mathbb{R}$ is the time-varying value of the $i$th agent, living on a one-dimensional manifold parametrized by $\mathbb{R}$, $\omega_i\in\mathbb{R}$ is the natural driving term of agent $i$, and $\xi_i$ will be used as an input to the system. 
Two agents $i$ and $j$ are interacting if a link between them exist in the interaction network, i.e., if and only if the corresponding term of the adjacency matrix $a_{ij}=1$. 
The interaction function between $i$ and $j$ is an odd, differentiable function $f_{ij}\colon\mathbb{R}\to\mathbb{R}$, and the coupling is symmetric, i.e., $f_{ij}=f_{ji}$. 
We also consider an attractive coupling, i.e., $\partial f_{ij}/\partial y > 0$ in an interval around $y=0$. 

If a fixed point $\bm{y}^*\in {\cal M}^n$ exists, one can linearize Eq.~\eqref{eq:dyn} around it, which yields, for a small deviation $\bm{x}=\bm{y}-\bm{y}^*$, to approximate the dynamics
\begin{align}\label{eq:dyn_lin}
 \dot{\bm{x}} &= -\J_{\bm f}(\bm{y}^*)\bm{x} + \bm{\xi} \, ,
\end{align}
where we use the Jacobian matrix of Eq.~\eqref{eq:dyn}, 
\begin{align}\label{eq:jac}
 \J_{\bm{f},ij}(\bm{y}^*) &= a_{ij}\frac{\partial}{\partial y}f_{ij}(y_i^*-y_j^*)\, .
\end{align}

One can verify that the oddness of the interaction ($f_{ij}$ odd) and the symmetry of the dynamics ($f_{ij}=f_{ji}$) implies that the Jacobian $\J_{\bm f}$ is a weighted Laplacian matrix of the interaction graph. 
It is then real symmetric, which implies that it has real eigenvalues, $\lambda_1\leq...\leq\lambda_n$, and its eigenvectors, $\bm{u}_1,...,\bm{u}_n$, form an orthonormal basis of $\mathbb{R}^n$. 
From now on, we will focus on stable fixed points of Eq.~\eqref{eq:dyn}, implying that the eigenvalues are nonnegative. 
Also, diffusive couplings imply that one eigenvalue vanishes, $\lambda_1=0$, with associated eigenvector $\bm{u}_1=n^{-1/2}(1,...,1)^\top$. 

{\bf Remark.}
{\it For other types of couplings where $\lambda_1>0$, the inference method proposed in Sec.~\ref{sec:net_inf} is even simpler than our case, we do not detail it. 
However, the result of Sec.~\ref{sec:n_inf} requires the zero mode and cannot be extended straighforwardly to nondiffusive couplings. }

Equation~\eqref{eq:dyn_lin} is then solved by expanding the deviation $\bm{x}$ over the eigenvectors $\bm{u}_\alpha$ of $\J_{\bm f}$, i.e., $x_i(t) = \sum_\alpha c_\alpha(t)u_{\alpha,i}$, with
\begin{align}\label{eq:sols}
 c_\alpha(t) &= e^{-\lambda_\alpha t}\int_0^te^{\lambda_\alpha t'}\bm{u}_\alpha^\top \bm{\xi}(t'){\rm d}t'\, .
\end{align} 
Details can be found, e.g., in~\cite{Tyl18a}. 

In order to get numerical confirmation of our results, we will apply them to the Kuramoto model on three different interaction graphs. 
The first one is a representation of the UK electrical grid (see inset in Fig.~\ref{fig:sine_ukersw}) with $n=120$ vertices and $m=165$ edges. 
The second one is a realization of an Erd\"os-Renyi graph with $n=120$ vertices and $m=329$ edges. 
The third one is a small-world graph realized according to the Watts-Strogatz process [\cite{Wat98}], with $n=120$ vertices and $m=242$ edges. 

\begin{figure*}
    \centering
    \includegraphics[width=.95\textwidth]{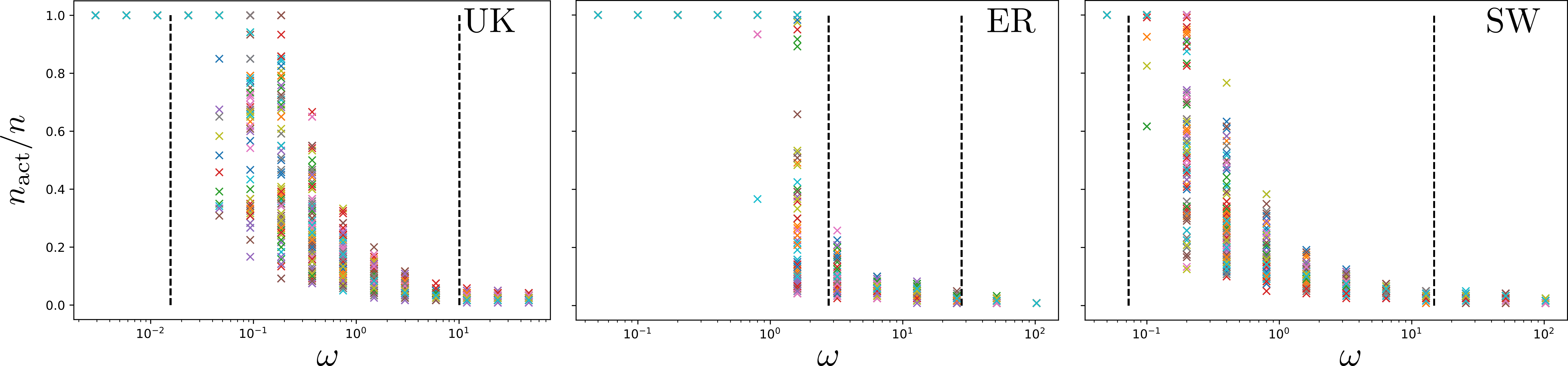}
    \caption{Fraction of responding agents vs. excitation frequency $\omega$ for the three networks UK (left panel), Erd\"os-Renyi (middle panel), and Small-World (right panel). 
    Vertical dashed lines correspond to $\lambda_2$ and $\lambda_n$\,.}
    \label{fig:spectra}
\end{figure*}

\section{Probing}
In order to determine oscillation modes in large electrical networks, one method is to apply a probing signals at some points of the network and measure the system's response at other points [e.g., \cite{Pie10}]. 
The probing is typically a sinusoidal signal with controlled amplitude and frequency. 

In the same spirit, we propose here to inject a sinusoidal signal at agent $i$ and to measure its impact at agent $j$. 
Let 
\begin{align}\label{eq:probing}
 \xi_i(t) &= a_0\sin(\omega_0t)\, ,
\end{align}
be the probing signal at agent $i$. 
We do not inject a probing signal at other nodes. 
Note that any signal shape could do the job, but the advantage of sine is that its amplitude and frequency are easily identifiable. 

To guarantee a minimal impact on the operation of the system, we need the amplitude $a_0$ to be sufficiently small. 
However, in most applications, the system under investigation will be subject to noise, and in order to be detectable, the probing amplitude should not be too small neither. 
Despite these contradicting requirements, we argue that, under our assumptions, an appropriate amplitude, satisfying both requirement simultaneously, can be chosen. 
Indeed, even though the system is subject to noise, we assume it is close to a steady state. 
If the noise amplitude was too large, it would push the system far away from its fixed point and it could not be considered as in (or close to) steady state. 
Under our assumptions, there is then a margin between the noise amplitude and the size of the basin of attraction of the steady state. 
The probing amplitude needs to be chosen in this margin. 
Determining the amplitude of the noise is rather straightforward from measurements. 
However, estimating the maximal tolerable disturbance magnitude preserving stability might represent a challenge and is beyong the scope of this manuscript [\cite{Men13}]. 

Also, the probing frequency $\omega_0$ needs to be kept small.
Keeping a small probing frequency guarantees that the system can adapt to the input and follow the probing signal. 
More precisely, a probing frequency can be qualified as small as long as it is smaller than the smallest eigenvalue of the Jacobian matrix ${\cal J}_{\bm f}$ in absolute value.

Introducing Eq.~\eqref{eq:probing} into Eq.~\eqref{eq:sols}, and recombining the eigenmodes yields the following response measured at agent $j$, 
\begin{align}\label{eq:recombine}
\begin{split}
 x_j^i(t) &=\sum_\alpha \frac{u_{\alpha,i}u_{\alpha,j} a_0}{\lambda_\alpha^2+\omega_0^2}  \\
 &\times\left[\lambda_\alpha\sin(\omega_0 t) + \omega_0 e^{-\lambda_\alpha t}- \omega_0 \cos(\omega_0 t) \right] \, .
\end{split}
\end{align}
To explicitly obtain an expression involving the Jacobian matrix, one should consider the long time limit $\lambda_\alpha t\gg 1$ with the asymptotic $\omega_0\ll\lambda_\alpha$ that yields,
\begin{align}\label{eq:trajs}
\begin{split}
 x_j^i(t) &= \left[\sum_{\alpha\ge 2}\frac{u_{\alpha,i} u_{\alpha,j}}{\lambda_\alpha}\right] a_0 \sin{(\omega_0 t)} + \frac{u_{1,i}^2 a_0}{\omega_0} [1-\cos(\omega_0 t)] \\
 &= \J_{\bm f,ij}^{\dagger}a_0 \sin{(\omega_0 t)} + \frac{a_0}{n\omega_0}[1-\cos{(\omega_0 t)]} \, ,
 \end{split}
\end{align}
where the $\dagger$ denotes the Moore-Penrose pseudo-inverse. 
We will say that a pair of agents $(i,j)$ can be probed if we have the ability to inject a probing signal at one of them and to measure the response at the other. 

Note that even though one usually has access to the dynamical variable $y_i$ instead of $x_i$, the latter can be recovered by comparing the behavior of $y_i$ before and after the introduction of the probing, 
\begin{align}
 x_i(t) &= y_i(t) - y_i^*\, .
\end{align}

\section{Spectrum range estimation}\label{sec:spec_range}
One condition in order to derive Eq.~\eqref{eq:trajs} is that the probing frequency $\omega_0$ is much smaller than all eigenvalues of the Jacobian matrix. 
From a practical point of view, however, the Jacobian being unknown, one cannot guarantee a priori to choose a frequency sufficiently small. 

Using the fact that we can probe and measure the system, we propose here a method to infer the range of the spectrum of the Jacobian Eq.~(\ref{eq:jac}). 
It can be inferred by tuning the signal's frequency $\omega_0$ and observing how many agents, $n_{\rm act}$, responded to the signal. 
Indeed, when $\omega_0\gg\lambda_\alpha$, one verifies in Eq.~\eqref{eq:recombine} that $x_j^i=O(\omega_0^{-1})$, which vanishes for large $\omega_0$, meaning that the signal stays local and do not spread across the network.
For $\omega_0\ll\lambda_\alpha$, Eq.~\eqref{eq:trajs} shows that for sufficiently small $\omega_0$, 
\begin{align}
 x_j^i(t) &\approx \frac{a_0}{n\omega_0}[1-\cos(\omega_0 t)]\, ,
\end{align}
independently of $i$ and $j$. 
All agents in the network then respond together to the signal (i.e., adiabatic shift in the parameters). 

Injecting a probing signal and varying its frequency, one can then estimate the range of the Jacobian's spectrum, which is delimited by the frequencies where probing propagates to the whole network and where it does not propagate at all respectively. 
This is illustrated in Fig.~\ref{fig:spectra}, where each cross corresponds to the fraction of agents that responded to the probing at a single agent. 
One observes that the boundaries (vertical dashed lines) of the spectrum of ${\cal J}_{\bm f}$ are roughly estimated by the transition from the probing signal staying local (large $\omega_0$) to spreading over the whole network (small $\omega_0$).

This gives at least an order of magnitude of the spectrum range, and hence of what is a "small" probing frequency.

\section{Number of agents}\label{sec:n_inf}
\begin{figure*}
    \centering
    \includegraphics[width=.95\textwidth]{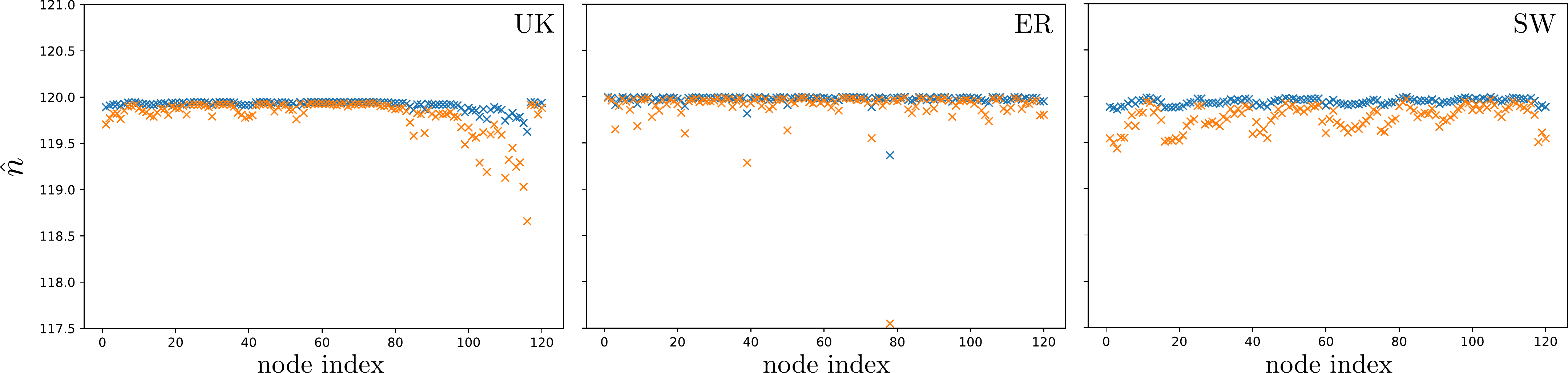}
    \caption{Estimation of the number of agents from a single probing/measurement for the three networks UK, Erd\"os-Renyi, and Small-World. 
    All of them have $n=120$ vertices. 
    Left panel, UK: $\omega_0/\lambda_2=0.02$ (orange), $\omega_0/\lambda_2=0.01$ (blue); 
    Middle panel, ER: $\omega_0/\lambda_2=0.0025$ (orange), $\omega_0/\lambda_2=0.00126$ (blue); 
    Right panel, SW: $\omega_0/\lambda_2=0.01$ (orange), $\omega_0/\lambda_2=0.005$ (blue)\,. 
    Some crosses are slightly higher than $120$ because of constant shifts in trajectories, as described in Sec.~\ref{sec:net_inf}. }
    \label{fig:node_number}
\end{figure*}

The number of agents in a coupled system is one of its primal properties. 
However, there are many physical examples where this number is not known exactly~[\cite{Su12}, \cite{Hae19}]
In this section, we introduce a method that allows to accurately determine the number of agents in any system governed by the dynamics of Eq.~\eqref{eq:dyn_lin}. 
Moreover, it requires to probe and measure a single node, which makes the method very efficient. 
Let us inject a probing signal, Eq.~\eqref{eq:probing} at node $i$, with $\omega_0\ll \lambda_\alpha$. 
Then for $\omega_0$ sufficiently small one has,
\begin{align}
 x_i^{\rm max}=\max_t|x_i(t)| &\approx \frac{2a_0}{n\omega_0} \, ,
\end{align}
from which one obtains an estimate for the number of nodes as,
\begin{align}\label{eq:node_est}
 \hat{n}=\frac{2a_0}{x_i^{\rm max}\omega_0} \, .
\end{align}
Note that we take the maximum to have a better accuracy in the estimation. However one can choose a particular time step $t$, keeping in mind that $t$ too short leads to	 vanishing values for Eq.~(\ref{eq:trajs}).
We check the validity of the estimation of Eq.~(\ref{eq:node_est}) in Fig.~\ref{fig:node_number}. Each cross corresponds to $\hat{n}$ obtained from a single node probing and measurement. One clearly sees that for $\omega_0$ small enough compared to $\lambda_2$, the estimated number of agents precisely matches the real one. 

Since the initial submission of this manuscript, this result has been developed and thoroughly discussed in~\cite{Tyl20}.

\section{Network inference}\label{sec:net_inf}
Assume now that each of the $[n(n-1)]/2$ pairs $(i,j)$, $1\leq i<j\leq n$, can be probed. 
For each pair $(i,j)$, we can measure either $\hat{x}_j^i(t)$ or $\hat{x}_i^j(t)$. 
According to Eq.~\eqref{eq:trajs} and by symmetry of ${\cal J}_{\bm f}$, these two trajectories should be the same, at least for $t$ sufficiently large, in order for the influence of initial conditions to vanish.

We will use the measured value $\hat{x}_j^i(t)$ at large $t$ to estimate ${\cal J}_{\bm f}^\dagger$ via Eq.~\eqref{eq:trajs}. 
To do this, we need to remember that the dynamics of Eq.~\eqref{eq:dyn} are invariant under a constant shift of all variables. 
This implies that two trajectories whose initial conditions differ only by a constant shift of all variables are exactly parallel, i.e., the constant shift is preserved for all time. 
This means that, depending on the initial conditions (unknown a priori), the measured trajectory $\hat{x}_j^i(t)$ might be shifted with respect to the predicted trajectory $x_j^i(t)$ of Eq.~\eqref{eq:trajs}. 
Denoting this constant shift as $c\coloneqq \hat{x}_j^i(t) - x_j^i(t)$, one can rearrange Eq.~\eqref{eq:trajs} as 
\begin{align}\label{eq:estim}
\begin{split}
 \hat{\cal J}^\dagger_{\bm{f},ij} &\coloneqq {\cal J}^\dagger_{\bm{f},ij} + \frac{c}{a_0\sin(\omega_0 t)} \\
 &= \left\{\hat{x}_j^i(t) - \frac{a_0}{n\omega_0}\left[1 - \cos(\omega_0 t)\right]\right\}/\left[a_0\sin(\omega_0 t)\right]\, .
\end{split}
\end{align}

Whereas we cannot determine the value of $c$, this is not an issue, as we show now. 
Indeed, the constant vector $\bm{u}_1=n^{-1/2}(1,...,1)^\top$ is an eigenvector of ${\cal J}_{\bm f}$, and then an eigenvector of ${\cal J}_{\bm f}^\dagger$ as well. 
As the eigenbasis of ${\cal J}_{\bm f}$ (and ${\cal J}_{\bm f}^\dagger$) is orthonormal, adding a constant value to each componenent of ${\cal J}_{\bm f}^\dagger$ does not modifies its eigenbasis, and modifies only one of its eigenvalues, $\lambda_1$. 
Diagonalizing $\hat{\cal J}_{\bm f}^\dagger$, it is then straightforward to replace its eigenvalue associated to $\bm{u}_1$ by $0$ (as it should be for the exact ${\cal J}_{\bm f}$ and ${\cal J}^\dagger_{\bm f}$) and to invert all other eigenvalues to recover ${\cal J}_{\bm f}$. 
Remark that in order to avoid singularity in Eq.~\eqref{eq:estim}, one should choose a time step $t$ such that $\sin(\omega_0 t)$ is sufficiently different from zero.

\begin{figure*}
    \centering
    \includegraphics[width=.95\textwidth]{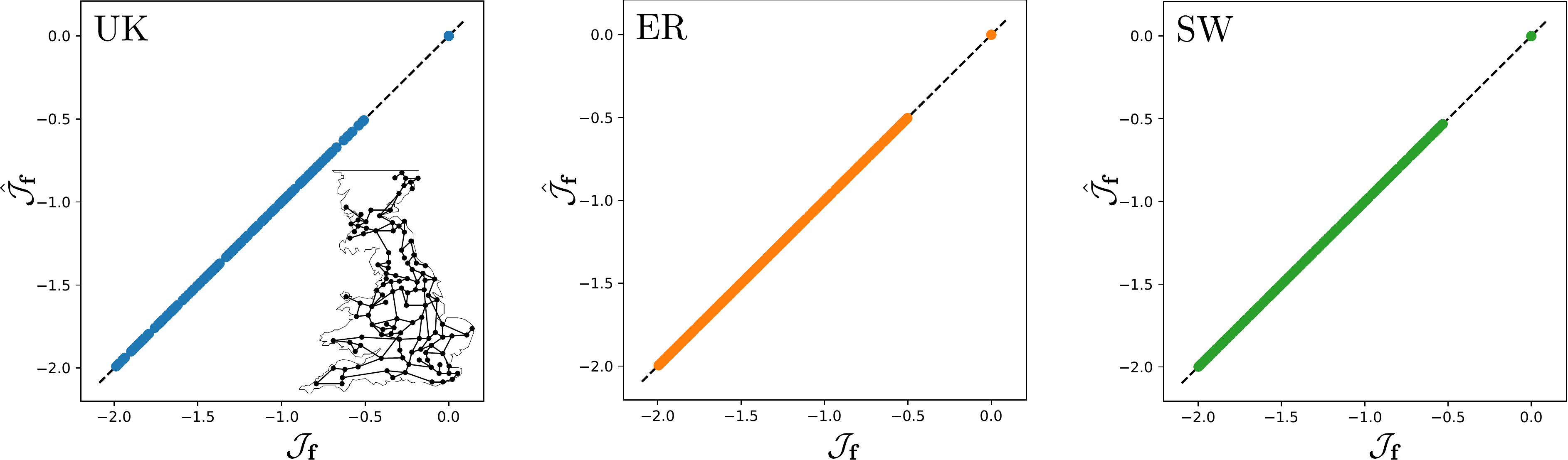}
    \caption{Inferred vs. real values of the elements of the Jacobian matrix ${\cal J}_{\bm f}$ for the three networks UK (left panel), Erd\"os-Renyi (middle panel), and Small-World (right panel).}
    \label{fig:sine_ukersw}
\end{figure*}

\begin{figure}
    \centering
    \includegraphics[width=\columnwidth]{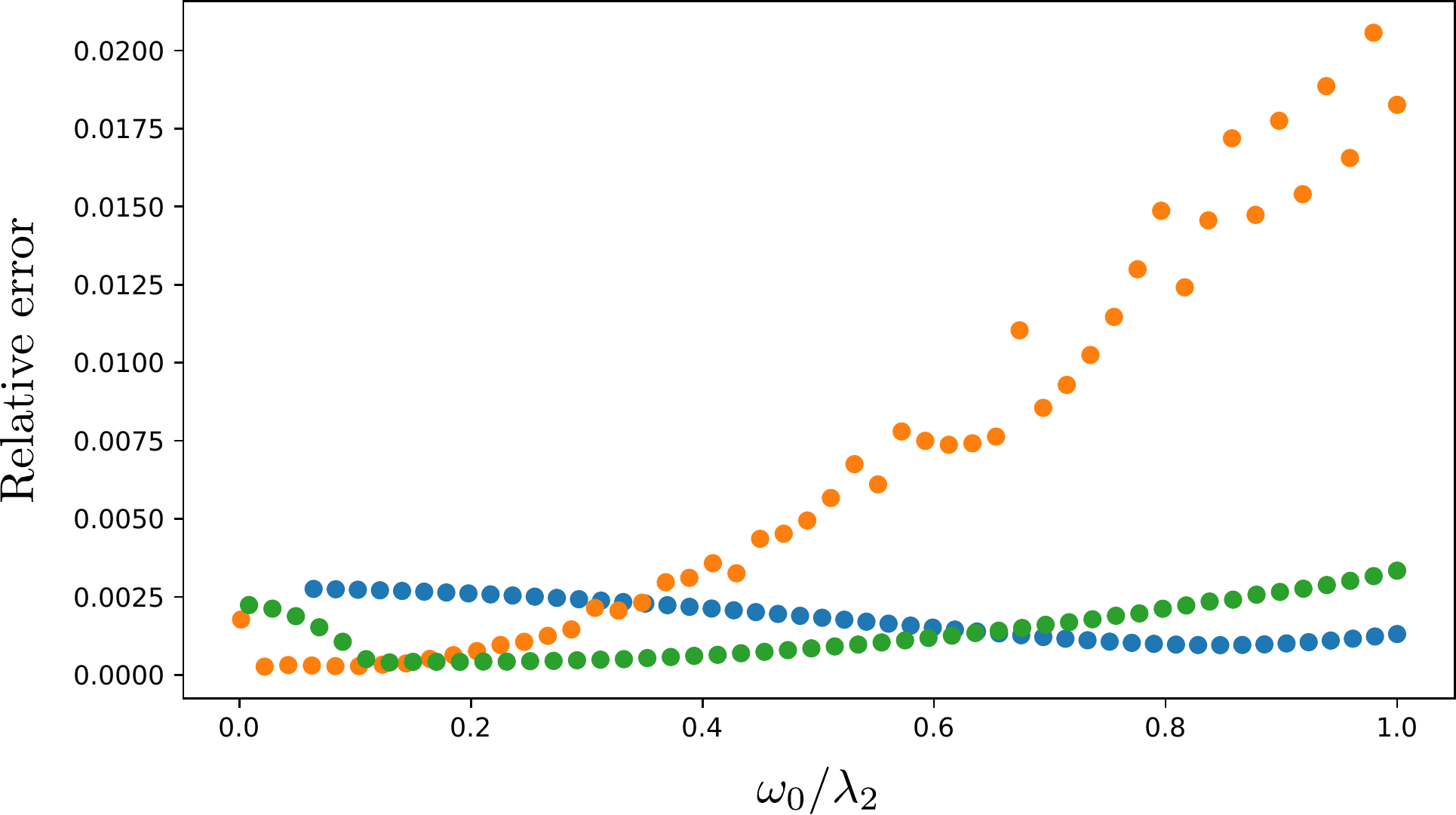}
    \caption{Relative error in the estimate of the Jacobian matrix $\hat{\cal J}_{\bm f}$ with respect to the frequency of the probing signal, normalized by the smallest eigenvalue of the Jacobian, for the three networks UK (blue), Erd\"os-Renyi (orange), and Small-World (green). 
    The relative error is computed as the Frobenius norm of the error of the Jacobian matrix normalized by the Frobenius norm of the real matrix. 
    One sees that for probing frequencies smaller that $10\%$ of $\lambda_2$, the estimate is very accurate. }
    \label{fig:nJ_vs_w0}
\end{figure}

The outcome of the procedure proposed above is illustrated in Fig.~\ref{fig:sine_ukersw}. 
Each dot corresponds to one element of the $120\times 120$ Jacobian matrix for the Kuramoto model on our selected graphs. 
One can see that the matrix is reconstructed with very high accuracy. 
In Fig.~\ref{fig:nJ_vs_w0}, we show the accuracy of the Jacobian estimate with respect to the frequency of the probing, normalized by $\lambda_2$. 
The accuracy of the estimate is measured as the Frobenius norm of the difference between the real Jacobian matrix and its estimate. 
We see that the accuracy is very good, even for probing frequencies close to $\lambda_2$. 

\section{Outlook}
We showed that based on a sinusoidal probing signal, injected at a node of a networked dynamical system while measuring the response of the network allows to:
\begin{itemize}
    \item Estimate the range of the spectrum of the Jacobian matrix of the system around its current fixed state;
    \item Estimate the number of agents in the system very efficiently: single node measurement; 
    \item Recover the network structure of the system with high fidelity.
\end{itemize}

The main advantages of our method are that it is non-invasive (does not disturb the system far from its operating state), it applies to a large set of coupling functions, whereas it requires to probe the system for a sufficiently long time, it only needs measurement at one time step.
Its main drawback is that, in order to recover the network using sinusoidal probing, one needs to be able to probe any pair of nodes in the network. 
The subsequent work would then investigate to what extent and with what confidence the network structure can be inferred if one has only probed a subset of the nodes.

\end{document}